\magnification=1200
\overfullrule=0pt
\centerline {{\bf A further improvement of a minimax theorem
of Borenshtein and Shul'man}}\par
\bigskip
\bigskip
\centerline {BIAGIO RICCERI}\par
\bigskip
\bigskip
\noindent
{\centerline {{\it Dedicated to Professor Ky Fan, with my greatest
admiration and esteem}}
\bigskip
\bigskip
\bigskip
\bigskip
The classical minimax problem can be formulated in the following
terms: given two non-empty sets $X, I$ and a function $f:X\times
I\to {\bf R}$, find suitable conditions ensuring the validity of
the equality
$$\sup_{\lambda\in I}\inf_{x\in X}f(x,\lambda)=\inf_{x\in X}
\sup_{\lambda\in I}f(x,\lambda)
\ . \eqno{(1)}$$
We refer to [2] and [6] for an overview of the subject.\par
\smallskip
 In the present paper, our starting point is the following result by
O. Yu. Borenshtein and V. S. Shul'man:\par
\medskip
THEOREM A ([1], Theorem 1). - {\it Let $X$ be a compact metric space,
$I$ a real interval, and $f:X\times I\to {\bf R}$ a continuous function
satisfying the following conditions:\par
\noindent
($a_{1}$)\hskip 5pt for every $x\in X$, the function $f(x,\cdot)$ is concave;
\par
\noindent
($a_{2}$)\hskip 5pt for every $\lambda\in I$, each local minimum of
the function $f(\cdot,\lambda)$ is a global minimum.\par
Then, equality $(1)$ holds.}\par
\medskip
Very recently, J. Saint Raymond obtained the  following remarkable
improvement of Theorem A:\par
\medskip
THEOREM B ([5], Theorem 1). - {\it Let $X, I$ be as in Theorem A, and let
$f:X\times I\to {\bf R}$ be a function satisfying the following conditions:
\par
\noindent
($b_{1}$)\hskip 5pt for every $x\in X$, the function $f(x,\cdot)$ is
concave and continuous;\par
\noindent
($b_{2}$)\hskip 5pt for every $\lambda\in I$, the function
$f(\cdot,\lambda)$ is
lower semicontinuous and each of its local minima is a global minimum.\par
Then, equality $(1)$ holds.}\par
\medskip
The aim of the present paper is to establish a further improvement
of Theorem A which improves Theorem B too. Our method of proof is
completely different from the ones of [1] and [5].\par
\smallskip
If $(X,\tau)$ is a topological space, we denote by $\tau_{s}$ the topology
on $X$ whose members are the sets $A\subseteq X$ such that $X\setminus
A$ is sequentially $\tau$-closed.
Clearly, $\tau_{s}$ is stronger than $\tau$. We note that a
function $f:X\to {\bf R}$ is sequentially $\tau$-lower semicontinuous
if and only
if it is $\tau_{s}$-lower semicontinuous. We also recall that a
real function $\varphi$ on a convex subset $C$ of a vector space is said to
be quasi-concave if,
for every $r\in {\bf R}$, the set $\{x\in C : \varphi(x)>r\}$ is convex.\par
\smallskip
Our result is as follows:\par
\medskip
THEOREM 1. - {\it Let $(X,\tau)$ be a topological space,
$\tau_{1}$ a topology on $X$ weaker than $\tau_{s}$, $I$ a real interval,
and $f:X\times I\to {\bf R}$ a function satisfying the following conditions:
\par
\noindent
($c_{1}$)\hskip 5pt for every $x\in X$, the function $f(x,\cdot)$
is quasi-concave and continuous;\par
\noindent
($c_{2}$)\hskip 5pt 
for every $\lambda\in I$, the function
$f(\cdot,\lambda)$ is
$\tau_{1}$-lower semicontinuous and each of its $\tau_{1}$-local minima
is a global minimum;\par
\noindent
($c_{3}$)\hskip 5pt there exist $\rho>\sup_{\lambda\in I}
\inf_{x\in X}f(x,\lambda)$ and $\lambda_{0}\in I$ such that the set
$$\{x\in X : f(x,\lambda_{0})\leq \rho\}$$
is sequentially $\tau$-compact.\par
Then, equality $(1)$ holds.}
\medskip
We now establish some intermediate results from which Theorem 1 will
easily follow.\par
\smallskip
We first recall that, if $T, S$ are two topological spaces, a multifunction
$F:T\to 2^S$ is said to be lower semicontinuous at $t\in T$ if, for every
open set $\Omega\subseteq S$, the condition $t\in F^{-}(\Omega)$ implies
that $t\in \hbox {\rm int}(F^{-}(\Omega))$, where $F^{-}(\Omega)=
\{t\in T: F(t)\cap \Omega\neq \emptyset\}$.\par
\smallskip
  Denote by $NLS(F(t))$ the set of all points
$s\in F(t)$ for which there exist a neighbourhood $\Omega$ of $s$
in $S$
and a net $\{t_{\alpha}\}$ in $T\setminus F^{-}(\Omega)$ converging to $t$.
It is clear that $F$ is not
lower semicontinuous at $t$ if and only if the set
$NLS(F(t))$ is non-empty.\par
\smallskip
 If $\tau$ is the topology on $S$, we will also
write $NLS_{\tau}(F(t))$ to denote $NLS(F(t))$.
\smallskip
The concept of lower semicontinuity is certainly one of the most important
in set-valued analysis. Thus, very often, it is of great interest to know
whether specific multifunctions are lower semicontinuous. 
The next proposition, on the contrary, does highlight a particular situation
where it is important to know that a certain multifunction is not
lower semicontinuous.\par
\medskip
PROPOSITION 1. - {\it Let $X, \Lambda$ be two topological spaces and
$f:X\times \Lambda\to {\bf R}$ a function such that, for each
$x\in X$, the set $\{\lambda\in \Lambda : f(x,\lambda)\geq 0\}$ is
closed. For each $\lambda\in \Lambda$, put
$$F(\lambda)=\{x\in X : f(x,\lambda)\leq 0\}\ .$$
Assume that the multifunction $F$ is not lower semicontinuous at
$\lambda_{0}\in \Lambda$.\par
Then, each point $x_{0}\in NLS(F(\lambda_{0}))$ is a local minimum
of the function $f(\cdot,\lambda_{0})$ 
such that $f(x_{0},\lambda_{0})=0$.}\par
\smallskip
PROOF. Let $x_{0}\in NLS(F(\lambda_{0}))$. So, there are a
neighbourhood $\Omega$ of $x_{0}$ in $X$ and a net
 $\{\mu_{\alpha}\}_{\alpha\in D}$
 in $\Lambda\setminus F^{-}(\Omega)$ converging to $\lambda_{0}$. Thus,
we have $f(x_{0},\lambda_{0})\leq 0$ as well as $f(x,\mu_{\alpha})>0$
for all $x\in \Omega$, $\alpha\in D$. Since each set
 $\{\lambda\in \Lambda : f(x,\lambda)\geq 0\}$ is closed, we then have
 $$f(x_{0},\lambda_{0})=0\leq f(x,\lambda_{0})$$
for all $x\in \Omega$, and the proof is complete.\hfill $\bigtriangleup$\par
\medskip
The next result is crucial.\par
\medskip
THEOREM 2. - {\it Let $(X,\tau)$ be a topological space,
$\tau_{1}$ a topology on $X$ weaker than $\tau_{s}$,
$I$ a real interval, and
 $F:I\to 2^X$ a multifunction, with $\tau_{1}$-closed values,
such that, for
each $x\in X$, the set $I\setminus F^{-}(x)$ is a non-empty interval open in
$I$. Finally,
suppose that there is $\lambda_{0}\in I$ such that $F(\lambda_{0})$
is non-empty and sequentially $\tau$-compact.\par
Then, the multifunction $F$ is not $\tau_{1}$-lower semicontinuous
at a point $\lambda^{*}\in I$ such that $NLS_{\tau_{1}}(F(\lambda^{*}))
\cap F(\lambda_{0})\neq \emptyset$.}\par
\smallskip
PROOF. For each $x\in X$, put
$$\Phi(x)=I\setminus F^{-}(x)\ .$$
Note that
$$\Phi^{-}(\lambda)=X\setminus F(\lambda)$$
for all $\lambda\in I$. Since $\tau_{1}$ is weaker than $\tau_{s}$,
each set $X\setminus \Phi^{-}(\lambda)$ is sequentially $\tau$-closed.
Hence, by
Proposition 2 of [4], we can find a compact (non-degenerate) interval
$[a,b]\subseteq I$, with $\lambda_{0}\in [a,b]$, such that
$$\Phi(x)\cap [a,b]\neq \emptyset$$
for all $x\in X$. Therefore, each set $\Phi(x)\cap [a,b]$ is a non-empty
interval open in $[a,b]$. For each $x\in X$, put
$$\alpha(x)=\inf (\Phi(x)\cap [a,b])$$
and
$$\beta(x)=\sup (\Phi(x)\cap [a,b])\ .$$
Let $x_{0}\in X$ and $\epsilon\in ]0,\beta(x_{0})-\alpha(x_{0})[$. Observe
that the set
$\Phi^{-}(\alpha(x_{0})+\epsilon)\cap \Phi^{-}(\beta(x_{0})-\epsilon)$
is a $\tau_{1}$-neighbourhood of $x_{0}$, and for every point $x$ belonging
to it, one has
$$\alpha(x)<\alpha(x_{0})+\epsilon$$
as well as
$$\beta(x_{0})-\epsilon<\beta(x)\ .$$
This shows that $\alpha$ is $\tau_{1}$-upper semicontinuous and that $\beta$
is $\tau_{1}$-lower semicontinuous in $X$. Now, suppose that $\lambda_{0}\in
]a,b[$. Then, we have
$$F(\lambda_{0})=\alpha^{-1}([\lambda_{0},+\infty[)\cup
\beta^{-1}(]-\infty,\lambda_{0}])\ .$$
Thus, since $F(\lambda_{0})$ is non-empty, either
$\alpha^{-1}([\lambda_{0},+\infty[)$ or
$\beta^{-1}(]-\infty,\lambda_{0}])$ is non-empty. For instance, suppose that
$\beta^{-1}(]-\infty,\lambda_{0}])$ is non-empty. Since $\tau_{1}$ is
weaker than $\tau_{s}$, the function $\beta$ is $\tau_{s}$-lower
semicontinuous and so, as we observed above, sequentially $\tau$-lower
semicontinuous. Hence, as $F(\lambda_{0})$ is sequentially $\tau$-compact,
there is $x^{*}\in X$ such that $\beta(x^{*})=\inf_{X}\beta$. Since
$\beta(x^{*})\leq \lambda_{0}$, we have $\beta(x^{*})\in ]a,b[$. This
implies, in particular, that $\beta(x^{*})$ does not belong to
$\Phi(x^{*})\cap [a,b]$, since this set is open in $[a,b]$. Hence, we have
$x^{*}\in F(\beta(x^{*}))$. Fix $c\in ]\alpha(x^*),\beta(x^*)[$.
Owing to the $\tau_{1}$-upper semicontinuity of $\alpha$,
the set $\alpha^{-1}(]-\infty,c[)$ is a $\tau_{1}$-neighbourhood of
$x^*$. Now, let $x\in \alpha^{-1}(]-\infty,c[)$ and $\lambda\in
]c,\beta(x^*)[$. Then, we have
$$\alpha(x)<c<\lambda<\beta(x^*)\leq \beta(x)\ ,$$
and so $\lambda\in \Phi(x)$, that is $x$ does not belong to $F(\lambda)$.
In other words, we have
$$F(\lambda)\cap \alpha^{-1}(]-\infty,c[)=\emptyset$$
for all $\lambda\in ]c,\beta(x^*)[$. This shows that
$x^{*}\in NLS_{\tau_{1}}(F(\beta(x^{*})))\cap F(\lambda_{0}).$
 In the case when
$\alpha^{-1}([\lambda_{0},+\infty[)$ is non-empty, $\alpha$ has at least one
global maximum, say $x'$. Proceeding as before (with
obvious changes), one realizes that
$x'\in NLS_{\tau_{1}}(F(\alpha(x')))\cap F(\lambda_{0}).$
 Suppose now that $\lambda_{0}=b$.
If $\inf_{X}\beta<b$, since
$$\beta^{-1}(]-\infty,b[)\subseteq F(b)$$
and $F(b)$ is sequentially $\tau$-compact,
$\inf_{X}\beta$ is attained by $\beta$, and, as we have seen, each
global minimum of $\beta$ belongs to
$NLS_{\tau_{1}}(F(\inf_{X}\beta))\cap F(\lambda_{0}).$
 If $\beta(x)=b$ for all
$x\in X$, pick $x_{0}\in F(b)$. Then, if $d\in ]\alpha(x_{0}),b[$ and
$\lambda\in ]d,b[$, the set $\alpha^{-1}(]-\infty,d[)$ is a
$\tau_{1}$-neighbourhood of $x_{0}$ which does not meet $F(\lambda)$.
So, $x_{0}\in NLS_{\tau_{1}}(F(b))$. In other words,
$F(b)=NLS_{\tau_{1}}(F(b))$.
Finally, if $\lambda_{0}=a$, reasoning as before (with obvious changes), one
checks that $\sup_{X}\alpha$
is attained by $\alpha$ and that each global maximum of $\alpha$ belongs to
$NLS_{\tau_{1}}(F(\sup_{X}\alpha))\cap F(a).$
 The proof is complete.\hfill $\bigtriangleup$\par
\medskip
REMARK 1. - From the above proof, it is clear that if we assume
$F(\lambda_{0})$ to be $\tau$-compact, then we get the same conclusion
provided that $\tau_{1}$ is weaker than $\tau$.\par
\medskip
A joint application of Proposition 1 and Theorem 2 gives the following
result:\par
\medskip
THEOREM 3. - {\it Let $(X,\tau)$ be a topological space, $\tau_{1}$ a
topology on $X$ weaker than $\tau_{s}$, $I$ a real interval and
$f:X\times I\to {\bf R}$ a function satisfying the following conditions:\par
\noindent
$(d_{1})$\hskip 5pt for each $x\in X$, the set $$\{\lambda\in I:
f(x,\lambda)\geq 0\}$$ is closed in $I$ and the set $$\{\lambda\in I:
f(x,\lambda)>0\}$$ is a non-empty interval open in $I$;\par
\noindent
$(d_{2})$\hskip 5pt for each $\lambda\in I$, the set
$$\{x\in X : f(x,\lambda)\leq 0\}$$ is $\tau_{1}$-closed; \par
\noindent
$(d_{3})$\hskip 5pt there is $\lambda_{0}\in I$ such that the set
$$\{x\in X : f(x,\lambda_{0})\leq 0\}$$
is non-empty and sequentially $\tau$-compact.\par
Then, there exists $\lambda^{*}\in I$ such that 
the function $f(\cdot,\lambda^{*})$ has a $\tau_{1}$-local minimum $x^{*}$
such that $f(x^{*},\lambda^{*})=0$ and $f(x^{*},\lambda_{0})\leq 0$.}\par
\smallskip
PROOF. For each $\lambda\in I$, put
$$F(\lambda)=\{x\in X : f(x,\lambda)\leq 0\}\ .$$
By Theorem 2,
the multifunction $F$ is not $\tau_{1}$-lower semicontinuous at a point
$\lambda^{*}\in I$ such that $NLS_{\tau_{1}}(F(\lambda^{*}))\cap
F(\lambda_{0})\neq \emptyset$. Any point of this set, by Proposition 1,
satisfies the conclusion.
 \hfill $\bigtriangleup$
\par
\medskip
REMARK 2. - Observe that Theorem 3 is a remarkable improvement of Theorem
1 of [3].\par
\medskip
{\it Proof of Theorem 1}. Arguing by contradiction, assume that
$$\sup_{\lambda\in I}\inf_{x\in X}f(x,\lambda)<\inf_{x\in X}
\sup_{\lambda\in I}f(x,\lambda)\ .$$
Then, if we fix $r$ satisfying
$$\sup_{\lambda\in I}\inf_{x\in X}f(x,\lambda)<r<
\min\left \{ \rho,\inf_{x\in X}
\sup_{\lambda\in I}f(x,\lambda)\right \} $$
and put
$$\varphi(x,\lambda)=f(x,\lambda)-r$$
for all $(x,\lambda)\in X\times I$, we see that the function $\varphi$
satisfies all the assumptions of Theorem 3. Consequently, there exist
$\lambda^{*}\in I$ and $x^{*}\in X$ such that $f(x^{*},\lambda^{*})=r$
and $x^*$ is a $\tau_{1}$-local minimum for $f(\cdot,\lambda^{*})$.
Since $\inf_{x\in X}f(x,\lambda^{*})<r$, $x^*$ is not a global minimum
for $f(\cdot,\lambda^{*})$, and this contradicts condition $(c_{2})$.
The proof is complete.\hfill $\bigtriangleup$\par
\bigskip
\bigskip
{\centerline {\bf References}}\par
\bigskip
\bigskip
\noindent
[1]\hskip 10pt O. YU. BORENSHTEIN and V. S. SHUL'MAN,
{\it A minimax theorem}, Math. Notes, {\bf 50} (1991), 752-754.\par
\smallskip
\noindent
[2]\hskip 10pt B. RICCERI and S. SIMONS (eds.), {\it Minimax theory
and applications}, Kluwer Academic Publishers, 1998.\par
\smallskip
\noindent
[3]\hskip 10pt B. RICCERI, {\it A new method for the study of nonlinear
eigenvalue problems}, C. R. Acad. Sci. Paris, S\'erie I, {\bf
328} (1999), 251-256. \par
\smallskip
\noindent
[4]\hskip 10pt B. RICCERI, {\it On a three critical points theorem},
 Arch. Math. (Basel), {\bf 75} (2000), 220-226.\par
\smallskip
\noindent
[5]\hskip 10pt J. SAINT RAYMOND, {\it On a minimax theorem}, Arch. Math.
(Basel), {\bf 74} (2000), 432-437.\par
\smallskip
\noindent
[6]\hskip 10pt S. SIMONS, {\it Minimax theorems and their proofs},
in Minimax and Applications, D.-Z. Du and P. M. Pardalos (eds.),
Kluwer Academic Publishers, 1995, 1-23.\par
\bigskip
\bigskip
Department of Mathematics\par
University of Catania\par
Viale A. Doria 6\par
95125 Catania, Italy\par
\smallskip
{\it e-mail address}: ricceri@dipmat.unict.it
\bye